\documentclass[11pt,twoside]{article}
\usepackage{amssymb}
\usepackage{mathrsfs}
\usepackage[plainpages=false]{hyperref}
\usepackage{indentfirst}
\usepackage{amscd}
\usepackage{amsfonts,latexsym,rawfonts,amsmath,amssymb,amsthm}
\usepackage{amsmath,amssymb,amsfonts,latexsym,lscape,rawfonts}

\newcommand{\Q}{\mathbb{Q}}
\newcommand{\R}{\mathbb{R}}
\newcommand{\C}{\mathbb{C}}
\newcommand{\N}{\mathbb{N}}
\newcommand{\Z}{\mathbb{Z}}

\newtheorem{prop}{Proposition}[section]
\newtheorem{thm}{Theorem}[section]
\newtheorem{lem}{Lemma}[section]
\newtheorem{claim}{Claim}

\newtheorem{rem}{Remark}[section]
\newtheorem{ex}{Example}[section]
\newtheorem{defi}{Definition}[section]

\newcommand{\beqs}{\begin{eqnarray*}}
\newcommand{\eeqs}{\end{eqnarray*}}

\begin{title}{On the
$\alpha$-Invariants of Cubic Surfaces with Eckardt Points}
\author {Yalong Shi}

\date{}

\end{title}

\begin{document}

\maketitle

\begin{abstract}
  In this paper, we show that the $\alpha_{m,2}$-invariant (introduced by Tian in \cite{[T4]} and \cite{[T5]}) of a smooth cubic
  surface with Eckardt points is strictly bigger than $\frac{2}{3}$. This can be used to simplify Tian's original proof of the
  existence of K\"ahler-Einstein metrics on such manifolds. We also
  sketch the computations on cubic surfaces with one ordinary double
  points, and outline the analytic difficulties to prove the existence of
  orbifold K\"ahler-Einstein metrics.
\end{abstract}


%\tableofcontents

\section{Introduction}

A very important problem in complex geometry is the existence of
canonical metrics, for example, the
 K\"ahler-Einstein metrics. An obvious necessary condition for the existence of
 a K\"ahler-Einstein
 metric is that the first Chern class of the manifold should be positive , zero or negative.
 Though the existence of K\"ahler-Einstein metrics was proved when $c_1 \leq 0$ by Aubin
(the $c_1<0$ case) and Yau (both the $c_1<0$ case and the $c_1=0$
case) in the 1970's, the $c_1>0$ case, i.e. the Fano case, is much
more complicated and still not completely understood today. However,
in complex dimension 2, the $c_1>0$ case is completely solved by
Tian in \cite{[T3]}. A complex surface with positive first Chern
class is also called a Del Pezzo surface. By the classification of
complex surfaces, the Del Pezzo surfaces are $\C P^2$, $\C P^1
\times \C P^1$, and $\C P^2$ blowing up at most eight points that
are in general positions. According to \cite{[T3]}, except for $\C
P^2$ blown up at one or two points, all the remaining Del
Pezzo surfaces admit K\"ahler-Einstein metrics.\\

Up to now, the most effective way to prove the existence of
K\"ahler-Einstein metrics is to use Tian's $\alpha$-invariants (and
$\alpha_G$-invariants for a compact group $G$) introduced in
\cite{[T1]}. We now recall the definitions.\\

Let $X$ be a compact K\"ahler manifold of dimension $n$ with
$c_1(X)>0$. $g$ is a K\"ahler metric with
$\omega_g:=\frac{\sqrt{-1}}{2\pi}g_{i\bar{j}}dz_i\wedge d \bar{z}_j
\in c_1(X)$. Define the space of K\"ahler potentials to be
$$P(X,g)=\{ \varphi\in C^2(X;\R) \big{|}\ \omega_g+\frac{\sqrt{-1}}{2\pi} \partial \bar{\partial} \varphi >0,\ \sup_{X} \varphi =0  \}.$$
We also define {\beqs
P_m(X,g) &=&\{\varphi\in P(X,g) \big{|}\ \exists\ a\ basis\ s_0,\dots,s_{N_m}\ of\ H^0(X,-mK_X),\\
&& s.t.\ \omega_g+\frac{\sqrt{-1}}{2\pi} \partial \bar{\partial}
\varphi = \frac{\sqrt{-1}}{2m\pi} \partial \bar{\partial}
 \log (|s_0|^2+\dots+|s_{N_m}|^2) \}.\eeqs}Here the functions $|s_i|^2$, $i=0, \dots, N_m$ are only defined
locally by choosing a local trivialization of $K_X^{-m}$. But it's
easy to see that the (1,1)-form on the right hand side is
independent of the trivialization we choose and is globally defined.

\begin{defi}[Tian \cite{[T1]}, \cite{[T2]}]
The $\alpha$-invariant and $\alpha_m$-invariant of $X$ are defined
to be:
$$\alpha(X):=\sup \{\alpha>0 \big{|}\  \exists\ C_\alpha>0,\ s.t.\ \int_X e^{-\alpha\varphi} dV_g \leq C_\alpha,\ \forall\ \varphi \in P(X,g)   \}$$
$$\alpha_m(X):=\sup \{\alpha>0 \big{|}\ \exists\ C_\alpha>0,\ s.t.\ \int_X e^{-\alpha\varphi} dV_g \leq C_\alpha,\
\forall\ \varphi \in P_m(X,g)  \}. $$
\end{defi}

If $G$ is a compact subgroup of $Aut(X)$, choose a $G$-invariant
K\"ahler form $\omega_g$, then we can define $P_G(X,g)$ and
$P_{G,m}(X,g)$ by requiring the potentials to be $G$-invariant.
Following the same procedure, we can also define the
$\alpha_G$-invariant and $\alpha_{G,m}$-invariant.\\

We have the following criteria of Tian:\footnote{Another very
convenient tool is the ``multiplier ideal sheaves" introduced by
Nadel in \cite{[N1]} and simplified by Demailly and Koll\'ar in
\cite{[DeKo]}. It's easy to see that their results are equivalent to
Tian's theorem, see \cite{[Shi]}.}

\begin{thm}[Tian \cite{[T1]}]
If $\alpha_G(X)>\frac{n}{n+1}$ where n is the complex dimension of
the Fano manifold X and G is a compact subgroup of Aut(X), then X
admits a G-invariant K\"ahler-Einstein metric.
\end{thm}

Recently, I. Cheltsov \cite{[Chel]} computed all the
$\alpha$-invariants and some of the $\alpha_G$-invariants for Del
Pezzo surfaces. Combined with Tian and Tian-Yau's earlier work, this
gives an alternative proof for the existence of K\"ahler-Einstein
metrics on all the Del Pezzo surfaces of degree less than 7 except
cubic surfaces with Eckardt points. Cheltsov showed that for cubic
surfaces with Eckardt points, the $\alpha$-invariants are exactly
$\frac{2}{3}$. One may ask whether we have $\alpha_G>\frac{2}{3}$
for some nontrivial group $G$ in the latter case. This is true for
some special cubic surfaces with Eckardt points, for example the
Fermat hypersurface in $\C P^3$. However, this is false in general.

\begin{ex}
 Let X be
the cubic surface defined by the equation
$$z_1^3+z_2^3+z_3^3+6z_1z_2z_3+z_0^2(z
_1+2z_2+3z_3)=0$$ where $[z_0,z_1,z_2,z_3]$ are the homogeneous
coordinates in $\C P^3$. Then according to \cite{[Dol]} Table 10.3,
$Aut(X)=\Z_2$, and $X$ has exactly one Eckardt point, namely
$[1,0,0,0]$. It's easy to see that the anti-canonical divisor cut
out by $z _1+2z_2+3z_3=0$ consists of the three coplanar lines and
is $\Z_2$ invariant. By the equivariant version of Theorem 2.2 (See
\cite{[De2]} for a proof.), we have $\alpha_{\Z_2}(X)=\frac{2}{3}$.
\end{ex}

So for all cubic surfaces with Eckardt points, the only known proof
for the existence of K\"ahler-Einstein metrics is still Tian's
original one in \cite{[T3]}. The key idea of Tian's proof in
\cite{[T3]} is to use his ``partial $C^0-$estimate". There are two
versions of ``partial $C^0-$estimate". The weaker one (Theorem 5.1
of \cite{[T3]}) states that the function
$$\psi_m=\frac{1}{m}\log (|s_0|^2_{\omega_{KE}}+\dots +|s_{N_m}|^2_{\omega_{KE}})$$
for any smooth K\"ahler-Einstein cubic surface $(X,\omega_{KE})$
satisfying $Ric(\omega_{KE})=\omega_{KE}$ has a uniform lower bound
for some $m$, where $\{s_i\}_{i=0}^{N_m}$ is an orthonormal basis of
$H^0(X,-mK_X)$ with respect to the inner product induced by
$\omega_{KE}$. The stronger one (Theorem 2.2 of \cite{[T3]}) says
that this holds for any sufficiently large $m$ satisfying $m\equiv
0\ (mod\ 6)$\footnote{This partial $C^0$-estimate should hold for
any sufficiently large m. See Conjecture 6.4 of \cite{[T6]}.}. If we
define the $\alpha_{m,2}$-invariant as follows:

\begin{defi}[Tian \cite{[T4]}, \cite{[T5]}]
Let $(X,\omega_g)$ be as above. The $\alpha_{m,2}$-invariant of $X$
is defined to be: \beqs \alpha_{m,2}(X) &:=& \sup\{\alpha>0 \big{|}\
\exists\ C_\alpha>0,\ s.t.\ \int_X
(|s_1|_g^2+|s_2|_g^2)^{-\frac{\alpha}{m}}
dV_g \leq C_\alpha, for\ any\\
&& s_1,s_2\in H^0(X,-mK_X)\ with\ \langle s_i, s_j
\rangle_g=\delta_{ij} \}. \eeqs
\end{defi}
Then Tian proved the following criterion:\\

\begin{thm}[Tian \cite{[T3]}, also\cite{[T4]}, \cite{[T5]}]
  Let $X$ be a smooth Del Pezzo surface obtained by blowing
  up $\C P^2$ at 5 to 8 points in general position. If for some integer $m\geq 0$, $\psi_m$ has a uniform lower
  bound on the deformations of $X$ that have K\"ahler-Einstein metrics,
  and
  $$\frac{1}{\alpha_m(X)}+\frac{1}{\alpha_{m,2}(X)}>3,$$
  then $X$ admits a K\"ahler-Einstein metric with positive scalar
  curvature.
\end{thm}

In the appendix of \cite{[T3]}, Tian proved that
$\alpha_{6k,2}>\frac{2}{3}$. Combining this with the stronger
version of ``partial $C^0-$estimate", he proved
the existence of K\"ahler-Einstein metrics on such manifolds using the above theorem.\\

In this paper, we prove the following theorem:
\begin{thm}[Main Theorem]
  Let X be a smooth cubic surface with Eckardt points, then for any
  integer $m>0$, we have $\alpha_{m,2}(X)>\frac{2}{3}$.
\end{thm}

One application of this theorem is to give a simplified proof of
Tian's theorem in \cite{[T3]}. With our theorem in hand, the weaker
version of ``partial $C^0-$estimate" is sufficient to prove the
existence of K\"ahler-Einstein metrics. We refer the reader
to \cite{[T4]} and \cite{[T5]} for more details.\\

The proof of the main theorem will be given in section 4. In section
2, we will discuss basic properties of Tian's invariants. Then we
compute the $\alpha$-invariant for cubic surfaces with Eckardt
points in section 3. This has already been done by Cheltsov in
\cite{[Chel]}. We include a direct proof here for the reader's
convenience. In section 5, we sketch the computations on cubic
surfaces with one ordinary double points since it is quite similar
to the smooth case. Then we outline an approach to establish the
existence of K\"ahler-Einstein metrics on singular stable cubic
surfaces, and discuss briefly the extra analytic difficulties in
this approach. The details will be presented elsewhere. We also
include an appendix on relations between $\alpha_m$-invariants and
the $\alpha$-invariant. This appendix is basically
taken from \cite{[Shi]}.\\

{\bf Acknowledgements:} The author would like to thank his thesis
advisor Professor Gang Tian for introducing him to this field and
for his continuous support. He thanks Professor Ivan Cheltsov for
pointing out a mistake in the first version of this paper and for
his useful comments on automorphism groups of cubic surfaces. He
thanks his teachers, Professor Weiyue Ding, Professor Xiaohua Zhu,
Professor Yuguang Shi and Professor Huijun Fan for their interest in
this work and for their encouragements. He also thanks his friends
Chi Li and Yanir Rubinstein for many helpful discussions. Part of
this work was done when the author was visiting Princeton University
in the academic year of 2007-2008. He would like to thank the
university and the Math Department for their hospitality, and thank
the Chinese Scholarship Council for their financial support. The
results of this paper
will be part of the author's Ph.D. thesis at Peking University.\\

\section{Preliminaries on Tian's invariants}

From the definitions of $\alpha$, $\alpha_m$ and $\alpha_{m,2}$, we
see that these invariants are not so easy to compute. In particular,
the uniform integration estimates involved in the definitions are
difficult to verify. Fortunately, we have the following
simi-continuity theorem for complex singularity exponents (= log
canonical thresholds):

\begin{thm}[Demailly-Koll\'ar \cite{[DeKo]}, also see Phong-Sturm \cite{[PhSt]}]
  Let X be a complex manifold. Let  P(X) be the set of locally $L^1$
  plurisubharmonic functions on X equipped with the topology of $L^1$
  convergence on compact subsets. Let K be a
  compact subset of X and define the complex singularity exponent by
  $$c_K(\varphi)=\sup \{ c>0 | e^{-2c\varphi}\in L^1_{loc}(U)\ for\ some
  \ neighborhood\ U\ of\ K\}.$$ If $\psi_i$ converges to $\varphi$ in
  P(X) and $c<c_K(\varphi)$, then $e^{-2c\psi_i}$ converges to $e^{-2c\varphi}$
  in $L^1(U)$ for some neighborhood U of K.
\end{thm}

The following proposition gives an alternative and easier way of
computing $\alpha_m(X)$:

\begin{prop} For any integer $m>0$, we have
$$\alpha_m(X)= \sup\{ c>0 | \int_X
|s|_{h^m}^{-\frac{2c}{m}} dV_g <+\infty, \forall\ s\in H^0(X,-mK_X),
s\neq 0\}.$$
\end{prop}

\begin{rem}
  If we define $c(s)$ to be the global complex singularity
exponent of $s$ (that is, $c(s)$ is the supremum of the set of
positive numbers c such that $|s|^{-2c}$ is globally integrable),
then the result of this proposition can be written as
$$\alpha_m(X)=\inf\{m\cdot c(s)| s\in H^0(X,-mK_X), s\neq 0 \}.$$ By
Theorem 2.1, we can actually find a holomorphic section $s\in
H^0(X,-mK_X)$ satisfying $\alpha_m(X)=m\cdot c(s).$\\
\end{rem}

\noindent {\bf Proof of Proposition 2.1:} We need only to show that
\begin{equation} \tag{$\clubsuit$}
\begin{split}
\alpha_m(X) &= \sup\{\alpha>0 \big{|}\ \exists\ C_\alpha>0,\ s.t.\
\int_X |s|_{h^m}^{-\frac{2\alpha}{m}} dV_g \leq C_\alpha, \\
& \quad \forall\ s\in H^0(X,-mK_X)\ with\ \int_X |s|^2_{h^m}dV_g=1
\}.
\end{split}
\end{equation}
 For if
this is true, the proposition follows easily from Theorem 2.1.\\

We now follow Tian's original
computations(\cite{[T2]}, \cite{[T3]}).\\

Assume the hermitian metric $h$ on $K^{-1}$ satisfies
$Ric(h)=\omega_g$. Fix an orthonormal basis $s_0,\dots, s_{N_m}$ of
$H^0(X,-mK_X)$ with respect to $h$ and $\omega_g$. For any $\varphi
\in P_m(X,g)$, there exists a basis of $H^0(X,-mK_X)$ $s'_0,\dots,
s'_{N_m}$, such that \beqs \omega_g+\frac{\sqrt{-1}}{2\pi} \partial
\bar{\partial} \varphi
&=& \frac{\sqrt{-1}}{2m\pi} \partial \bar{\partial} \log (|s'_0|^2+\dots+|s'_{N_m}|^2)\\
&=& \frac{\sqrt{-1}}{2m\pi} \partial \bar{\partial} \log
(|s'_0|_{h^m}^2+\dots+|s'_{N_m}|_{h^m}^2)+\omega_g \eeqs Set
$$\tilde{\varphi}:= \frac{1}{m} \log
(|s'_0|_{h^m}^2+\dots+|s'_{N_m}|_{h^m}^2),$$ then
$\varphi=\tilde{\varphi}-\sup_X \tilde{\varphi}$. Since the value of
$\tilde{\varphi}$ doesn't change under unitary transformations on
$H^0(X,-mK_X)$ with respect to $h$ and $\omega_g$, we may assume
further that $s'_i=\lambda_i s_i$, $0<\lambda_0 \leq \dots \leq
\lambda_{N_m}$. So we can write $\tilde{\varphi}$ as
$$\tilde{\varphi}= \frac{1}{m} \log
(\lambda_0^2|s_0|_{h^m}^2+\dots+\lambda_{N_m}^2|s_{N_m}|_{h^m}^2).$$
Observe that $$\sup_X
\tilde{\varphi}=\frac{1}{m}\log(\lambda_{N_m}^2)+O(1) ,$$ we can
write
$$\varphi=\frac{1}{m} \log (\sum_{i=0}^{N_m}\lambda_i^2|s_i|_{h^m}^2)+O(1)$$
with $0<\lambda_0 \leq \dots \leq \lambda_{N_m}=1$. Then the
equality $(\clubsuit)$ follows easily from this expression and
Theorem 2.1. $\Box$\\

\begin{thm}[Demailly \cite{[De2]}] For any Fano
manifold $X$, we have
 $$\alpha(X)=\inf\{\alpha_m(X)\}
=\lim_{m\to +\infty} \alpha_m(X).$$
\end{thm}

A simple proof of Theorem 2.2 will be given in the appendix, which
is basically taken from the author's thesis in preparation. Note
that though Demailly's proof looks more complicated than the proof
we give here, his proof can yield more information in the
equivariant case. We
refer the interested readers to his paper for more details. \\

\begin{rem}
  A conjecture of Tian claims that for any Fano manifold X,
  one has $\alpha(X)=\alpha_m(X)$ when
  m is sufficiently large.
  We will discuss this problem in a separate paper.
\end{rem}

Now we state a similar proposition for $\alpha_{m,2}$-invariants,
whose proof is quite easy and thus omitted.
\begin{prop}
Let $(X,\omega_g)$ be as in Definition 1.1. Then we have: \beqs
\alpha_{m,2}(X) &=& \sup\{\alpha>0 \big{|}\  \int_X
(|s_1|_g^2+|s_2|_g^2)^{-\frac{\alpha}{m}}
dV_g <+\infty, for\ any\\
&& s_1,s_2\in H^0(X,-mK_X)\ with\ \langle s_i, s_j
\rangle_g=\delta_{ij} \}.\\\eeqs
\end{prop}

\section{The $\alpha$ invariants of cubic surfaces with Eckardt points}

Let X be a smooth cubic surface in $\C P^3$. It's well known that
there are exactly 27 lines on X. If we realize X as $\C P^2$ blowing
up 6 generic points $p_1,\dots,p_6$, then the 27 lines are:
\begin{itemize}
  \item the exceptional divisors: $E_1, \dots, E_6$;
  \item the strict transforms of lines passing through 2 of the 6
  points: $L_{12}, \dots, L_{56}$;
  \item the strict transforms of the quadrics that avoids only 1 of
  the 6 points: $F_1, \dots, F_6$.
\end{itemize}

It's easy to check that each line above intersects with other 10
lines, and that if 2 lines intersect, then there is a unique other
line that intersects them both. If it happens that there are  three
coplanar lines intersecting at one point $p$ on X, then we call $p$
an ``Eckardt point". Note that a generic cubic surface does not have
any Eckardt points. For detailed information about cubic surfaces,
we refer
the reader to the books \cite{[GH]}, \cite{[Hart]} and \cite{[Dol]}.\\

We shall prove the following theorem of Cheltsov in this section.

\begin{thm}[Cheltsov\cite{[Chel]}]
 Let X be a smooth cubic surface with Eckardt points, then for any
 integer $m>0$, $\alpha_m(X)=\alpha_1(X)=\frac{2}{3}$. In
 particular, by Theorem 2.2, $\alpha(X)=\frac{2}{3}$.
\end{thm}

\subsection{Computing $\alpha_1(X)$ }

 To compute the $\alpha_1$-invariant of our cubic surface, by Proposition
 2.1, we need
only to consider the singularities cut out by anti-canonical
sections. This is done, for example, in \cite{[T3]} and
\cite{[Par]}. The most ``singular" sections are exactly those
defined by triples of lines intersecting at Eckardt points. It's
easy to see that the singularity exponents of these sections
are equal to $\frac{2}{3}$.\\

\subsection{Computing $\alpha(X)$ }

Now we show that $\alpha(X)=\alpha_1(X)$. The main tool is the
following theorem:\footnote{Actually, what we use in this paper is
just the connectness of ``multiplier ideal subschemes", which in
fact can also be proved directly by H\"ormander's $L^2$ method, see
\cite{[T4]}. }

\begin{thm}[Nadel vanishing theorem]
  Let X be a smooth complex projective variety, let D be any
  $\Q$-divisor on X, and let L be any integral divisor such that
  $L-D$
  is nef and big. Denote the multiplier ideal sheaf of D by $\mathcal
{J}(D)$, then $$H^i(X,\mathcal {O}_X(K_X+L)\otimes \mathcal
{J}(D))=0$$for any $i>0$.\\
\end{thm}

\noindent {\bf Proof of Theorem 3.1:} Suppose
$\alpha(X)<\alpha_1(X)$, then there is an integer $m$ such that
$\alpha_m(X)<\alpha_1(X)$. Then by definition, there is a nontrivial
holomorphic section $s$ of $K_X^{-m}$£¬such that
$c(s)<\frac{\lambda}{m}$, where $\lambda \in \Q$ and
$\lambda<\frac{2}{3}$. We denote the corresponding effective divisor
by $Z(s)\in
|-mK_X|$.\\

Now we need a lemma:

\begin{lem}
  For any $0<\lambda\leq \frac{2}{3}$ and any nonzero holomorphic section
  $s\in H^0(X,-mK_X)$, the locus of non-integrable points of
  $|s|^{-\frac{2\lambda}{m}}$ is a
  single point.
\end{lem}

\noindent {\bf Proof of Lemma 3.1:} Denote by $Z(s)$ the effective
divisor defined by the section $s$. Apply Nadel's theorem to the
sheaf $\mathcal {J}(\frac{\lambda}{m}Z(s))$ and the integral divisor
$-K_X$, we know that the locus of non-integrable points of
$|s|^{-\frac{2\lambda}{m}}$ should be a connected subset of $X$. We
denote the locus by $C$. So we need only to show that $C$ does not
contain one dimensional parts.\\

Suppose this is not true. Write $C=\cup_i  C_i$, where the $C_i$'s
are different irreducible curves. Then we can write $Z(s)$ as
$$Z(s)=\sum_i \mu_i C_i+ \Omega$$
where $\Omega$ is an effective divisor whose support doesn't contain
any of the $C_i$'s and the $\mu_i$'s are integers such that
$\frac{\lambda
\mu_i}{m}\geq 1$.\\

First, we know that every $C_i$ is a smooth rational curve, for if
not, we have $C_i^2>0$, then
$$0\leq 2p_a(C_i)-2=(C_i+K_X)\cdot C_i \leq (1-\frac{\mu_i}{m})C_i^2
-\sum_{j\neq i}\frac{\mu_j}{m}C_j\cdot C_i<0$$ A contradiction.\\

On the other hand, we have
$$3m=-Z(s)\cdot K_X=-\sum_i \mu_i C_i\cdot K_X-K_X \cdot \Omega\geq \frac{3m}{2}\sum_i C_i\cdot (-K_X)$$
So there are three possibilities:
\begin{enumerate}
  \item $\Omega$ is empty and there are two $C_i$'s, both among the 27 lines;
  \item $\Omega$ is empty and there is only one $C_i$, with $C_1^2=0$ and $K_X\cdot C_1=-2$;
  \item There is only one $C_i$,and it is one of the 27 lines.
\end{enumerate}

In case 1, $\lambda=\frac{2}{3}$, $\mu_1=\mu_2=\frac{3m}{2}$, and
$Z(s)=\frac{3m}{2}(C_1+C_2)$. But $(C_1+C_2)$ can not be an ample
divisor. A contradiction.\\

 In case 2, $\lambda=\frac{2}{3}$,
$\mu=\frac{3m}{2}$, and $Z(s)=\frac{3m}{2}C_1$. Since $C_1^2=0$,
$C_1$ can not be ample, a contradiction.\\

In case 3, write $Z(s)=\mu C_1+\Omega$. Choose a birational morphism
$\pi$ from X to $\C P^2$ such that $deg \pi (C_1)=2$. Then
$$3\lambda= \pi^{*} H \cdot \frac{\lambda}{m} Z(s)\geq \pi^{*} H \cdot
 \frac{\lambda}{m}\mu C_1=H\cdot \frac{\lambda}{m}\mu \pi(C_1)= \frac{2\lambda}{m}\mu\geq2 $$
If $\lambda<\frac{2}{3}$, this is already a contradiction. If
$\lambda=\frac{2}{3}$, then $\mu=\frac{3m}{2}$ and $\Omega$ consists
of exceptional divisors, i.e. lines. Write
$$Z(s)=\frac{3m}{2} C_1+\sum_i \kappa_i L_i.$$
If $L$ is a line intersects with $C_1$, then $L$ must be contained
in $\Omega$, for otherwise $m=L\cdot Z(s)\geq \frac{3m}{2}$. A
contradiction. So $\Omega$ must contain the 10 lines having positive
intersection numbers with $C_1$. On the other hand, if $L_i\cdot
C_1=1$, then
$$m=L_i\cdot Z(s)\geq \frac{3m}{2}-\kappa_i,$$
so $\kappa_i\geq \frac{m}{2}$. However,
$$m=C_1\cdot Z(s)=-\frac{3m}{2}+\sum_i \kappa_i L_i\cdot C_1$$
So there are at most 5 $L_i$'s having positive intersection numbers
with $C_1$. A contradiction.                                 $\Box$\\

\begin{rem}
  The above lemma actually holds for any $0<\lambda<1$. We refer the reader to \cite{[N2]} for a proof.
\end{rem}

Now we continue to prove Theorem 3.1. Denote the point in the above
lemma by $p$. Now choose a birational morphism $\pi$ from X to $\C
P^2$ such that it is an isomorphism near $p$. Then $\pi (Z(s))$ is
an effective divisor of $\C P^2$. It's obvious that $\pi (Z(s)) \in
|-mK_{\C P^2}|$. Choose a generic line $L$ of $\C P^2$ that doesn't
pass $\pi(p)$. Let's now consider the $\Q$-divisor
$\Omega:=\frac{\lambda}{m}\pi(Z(s))+L$ which is numerically
equivalent to $(3\lambda +1)H$. Consider the multiplier ideal sheaf
$\mathcal {J}(\Omega)$. By Nadel's vanishing theorem, the multiplier
ideal subscheme associated with $\mathcal {J}(\Omega)$ should be
connected. But from our construction, its support should be
$\{\pi(p)\} \cup L$, which is obviously not connected. A
contradiction.  $\Box$\\

\section{Proof of the main theorem }

The key to the proof of the main theorem is the following:

\begin{thm}
  Suppose X is a smooth cubic surface in $\C P^3$ with Eckardt
  points, m is a positive integer, then if $s\in H^0(X,-mK_X)$ is a section such that
  $c(s)=\frac{2}{3m}$, then there exists a section $s_1\in
  H^0(X,-K_X)$ such that $s=s_1^{\otimes m}$. Moreover, the support
  of $Z(s_1)$ consists of three lines intersecting at one of the
  Eckardt points.
\end{thm}

Actually, Tian proved the theorem in the case of $m=6$ in the
appendix of \cite{[T3]}. Our proof is greatly inspired by his. The
proof is based on the following observations:

\begin{lem}
  Suppose $s=s_1\otimes s_2$ with $s_1 \in H^0(X, -K_X)$, $s_2 \in H^0(X,
  -(m-1)K_X)$ and $c(s)=\frac{2}{3m}$, then $c(s_1)=\frac{2}{3}$ and
  $c(s_2)=\frac{2}{3(m-1)}$.
\end{lem}

\noindent {\bf Proof:}$\quad$ We have $c(s_1)\geq\frac{2}{3}$ and
  $c(s_2)\geq \frac{2}{3(m-1)}$ by Theorem 3.1. Then the lemma is trivial by
  the H\"older inequality.   $\Box$\\

\begin{lem}
If $s\in H^0(X,-mK_X)$ is a holomorphic section and $mult_p s > m$
for some point p. If p lies on one of the lines $L$, then $L$ is
contained in the support of $Z(s)$. In particular, if $mult_p s>m $,
and p is an Eckardt point, then the three lines through p are all
contained in the support of $Z(s)$, hence $s=s_1\otimes s_2$ with
$s_1 \in H^0(X, -K_X)$, $s_2 \in H^0(X,
  -(m-1)K_X)$.
\end{lem}

\noindent {\bf Proof:}$\quad$Suppose $L$ is not contained in the
support of $Z(s)$, then
$$m=L\cdot Z(s) \geq mult_p s >m$$
A contradiction. $\Box$\\

\begin{lem}
  If $f\in \mathcal{O}_{\C^2,0}$ with $mult_0 f=k$, then $\frac{1}{k} \leq c_0(f)\leq \frac{2}{k}
  $, and if moreover $c_0(f)=\frac{1}{k}$, then locally $f=gh^k$,
  where $g(0)\neq 0$ and $mult_0 h =1$.
\end{lem}

\noindent {\bf Proof:}$\quad$ This lemma follows easily from the
fact that in dimension two, one can compute the singularity exponent
via Newton polygons for some analytic coordinates. We refer the
reader to Varchenko's paper \cite{[Var]}, the appendix of Tian's
paper \cite{[T3]} and the book of Koll\'ar, Smith and Corti
\cite{[Kol-Sm-Co]} for
detailed proofs. $\Box$ \\

 Based on these lemmas, we need only to show
that the only point $p$ where $c_p(s)=\frac{2}{3m}$ is an Eckardt
point. Actually the arguments in \cite{[Chel]} already imply this,
but his proof is more complicated and uses some properties of Geiser
involutions. So we give a simple proof here,
 which avoids Geiser involutions but still uses some observations of
 Cheltsov \cite{[Chel]} and Tian \cite{[T3]}.\footnote{Theorem 4.1
 here is a special case of Theorem 4.1 of Cheltsov \cite{[Chel2]}, which was not
 known to me when I wrote the first version of this paper. Here I
 give an alternative and more direct proof. The readers can consult \cite{[Chel2]}
 for the proof of a more general result.}\\

\noindent {\bf Proof of Theorem 4.1:} Suppose $p$ is not an Eckardt
point, then there are three possibilities:
\begin{enumerate}
  \item $p$ doesn't belong to any of the lines;
  \item $p$ belongs to exactly one line;
  \item $p$ belongs to exactly two lines.
\end{enumerate}

We shall rule them out one by one. First, note that by Lemma 3.1 and
Lemma 4.3, we have
 $mult_p s >\frac{3m}{2}$.\\

{\bf Case 1}

 In this case, we may choose a $D\in |-K_X|$ which has multiplicity
 2 at $p$. Then if D is not contained in the support of $Z(s)$, we have
 $$3m=D\cdot Z(s) \geq 2 mult_p s > 3m.$$
 A contradiction.\\

{\bf Case 2}

 In this case, there is a section $s'\in H^0(X,-K_X)$
such that $Z(s')=L_1+D$ where $L_1$ is a line and $D$ is a quadratic
curve intersecting with $L_1$ at $p$. Then $L_1$ is contained in the
support of $Z(s)$, so $D$ is not contained in the support of $Z(s)$
in view of Lemma 4.1. Thus
$$2m=D\cdot Z(s)\geq mult_p D \cdot mult_p s > \frac{3m}{2}\cdot mult_p D,$$
which implies $mult_p D=1$.\\

Write $Z(s)=\mu L_1+\Omega$, where $\Omega$ is an effective divisor
whose support doesn't contain $L_1$. We have
$$m=L_1 \cdot Z(s) \geq -\mu+mult_p \Omega = mult_p s - 2\mu$$
and
$$2m=D\cdot Z(s) \geq 2\mu + mult_p \Omega = mult_p s + \mu,$$
which imply $mult_p s \leq \frac{5m}{3}$ and
$\frac{m}{4}<\mu<\frac{m}{2}$.\\

We blow up $X$ at $p$ to get a surface $U$, $\pi : U \to X$. For any
divisor $F$ of $X$, denote by $\bar{F}$ the strict transform of $F$.
We have
$$\pi^* (K_X+ \frac{2}{3m}Z(s))=K_U+ \frac{2}{3m}(\mu \bar{L}_1+\bar{\Omega})+(\frac{2}{3m}mult_p s-1)E.$$
Since the pair $(X,\frac{2}{3m}Z(s))$ is not log terminal at $p$,
there is a point $Q\in E$ satisfying $$\frac{2}{3m}(\mu \ mult_Q
\bar{L}_1+mult_Q \bar{\Omega})+(\frac{2}{3m}mult_p s-1)\geq 1,$$
that is
$$\mu \ mult_Q \bar{L}_1+mult_Q \bar{\Omega}+mult_p s \geq 3m.$$

Now we prove that $Q \notin \bar{L}_1$. If $Q\in \bar{L}_1$, then
since $\frac{2\mu}{3m}<1$, by the Fubini theorem and Theorem 2.1, we
know that $$\bar{L}_1 \cdot [\frac{2}{3m}
\bar{\Omega}+(\frac{2}{3m}mult_p s-1)E] \geq 1,$$ which implies that
$\mu>m$. A contradiction.\\

So we have $$mult_Q \bar{\Omega}+mult_p s \geq 3m.$$Then we can see
that $mult_Q \bar{\Omega}\geq \frac{4m}{3}$ and hence $mult_p \Omega
\geq \frac{4m}{3} $. Moreover, it's easy to see that actually
$mult_p
\Omega=\frac{4m}{3}$ and $\mu= \frac{m}{3}$.\\

If m is not a multiple of 3, this already leads to a contradiction.
If $m=3k$, then $Z(s)=kL_1+\Omega$ with $mult_p \Omega=4k$. In this
case, $L_1$ intersects with $\Omega$ only at $p$, and $L_1$ is not
tangent to $\Omega$ at $p$. By Theorem 2.1 and the Fubini theorem,
the singularity exponent of s at $p$ is at least $\frac{1}{4k}$
which is bigger than
$\frac{2}{3m}=\frac{2}{9k}$. A contradiction.\\

{\bf Case 3}

In this case, there is a section $s' \in H^0(X,-K_X)$, with
$Z(s')=L_1+L_2+L_3$ where $L_1$ and $L_2$ intersects at $p$ and
$L_3$ is the other line coplanar with $L_1, L_2$ and $p\notin L_3$.
Firstly $L_1$ and $L_2$ must be contained in the support of $Z(s)$
as
before.\\

By Lemma 4.1, $L_3$ is not contained in the support of $Z(s)$. Write
$Z(s)=\mu L_1 +\nu L_2 + D$, then
$$m=L_1\cdot Z(s) =-\mu+\nu +L_1\cdot D \geq -\mu+\nu +mult_p D$$
$$m=L_2\cdot Z(s) =\mu-\nu +L_2\cdot D\geq \mu-\nu +mult_p D$$
$$m=L_3\cdot Z(s) = \mu+\nu +L_3 \cdot D \geq \mu+\nu$$
These imply that $$\mu+\nu\leq m, \ \frac{m}{2}< mult_p D \leq
m\quad and\ \mu> \frac{m}{4},\ \nu>\frac{m}{4}.\\$$

As in Case 2, we blow up $X$ at $p$ to obtain a surface $U$. We have
$$\pi^* (K_X+\frac{2}{3m}Z(s))=K_U + \frac{2}{3m}
(\mu \bar{L}_1+\nu \bar{L}_2+ \bar{D})+(\frac{2}{3m}mult_p s -1)E.$$
As before, there is a point $Q$ on $E$ satisfying
$$\mu\ mult_Q \bar{L}_1 + \nu\ mult_Q \bar{L}_2 + mult_Q \bar{D} + mult_p s \geq 3m.$$
It's easy to see that $Q \notin \bar{L}_1 \cup \bar{L}_2$, so the
above inequality reduces to (with $\mu+\nu \leq m $ in mind)
$$mult_Q \bar{D} + mult_p D \geq 2m.$$
Since $mult_p D \leq m$, we must have $\mu=\nu=\frac{m}{2}$, and
$mult_p D=m$. \\

If $m$ is odd, this already leads to a contradiction. Now suppose
$m=2k$. We can write the section s locally as $s=z_1^k z_2^k h$,
with $mult_0 h = 2k$. By the H\"older inequality and the fact that
$$c_0(s)=\frac{1}{3k},\quad c_0(z_1^k z_2^k)=\frac{1}{k},$$we know that
$c_0(h)\leq \frac{1}{2k}$. So by Lemma 4.3, we can write $\Omega$
locally at $p$ as $\Omega= 2k C$, where $C$ is a curve regular at
$p$ and not tangent to $L_1$ or $L_2$. Then by blowing up $p$ we get
a log resolution for the pair $(X, Z(s))$ near $p$. It's easy to see
that $$c_p(s)=lct_p(X, Z(s))=\frac{1}{2k}>
\frac{1}{3k}=\frac{2}{3m}.$$A contradiction. $\Box$ \\

Now let's turn to the proof of the main theorem:\\

\noindent {\bf Proof of Theorem 1.2:} Obviously,
$\alpha_{m,2}(X)\geq \alpha_m(X)=\frac{2}{3}$. It's easy to show
that there are orthogonal sections $s_1,s_2 \in H^0(X,-mK_X)$ such
that
$$\alpha_{m,2}(X)=m c((|s_1|_g^2+|s_2|_g^2)^{\frac{1}{2}}).$$ To prove
the theorem, by compactness arguments, it suffices to show that at
every point $p\in X$,
$$c_p((|s_1|_g^2+|s_2|_g^2)^{\frac{1}{2}})>\frac{2}{3m}.$$

This is simple. Suppose
$$c_p((|s_1|_g^2+|s_2|_g^2)^{\frac{1}{2}})=\frac{2}{3m}$$ for some point $p$. By
comparing $(|s_1|_g^2+|s_2|_g^2)^{\frac{1}{2}}$ with $|s_1|_g$ and
$|s_2|_g$ respectively, we know that $c_p(s_i)=\frac{2}{3m}$. Hence
by Theorem 4.1, $p$ must be an Eckardt point and $s_1=\lambda
s_2$ for some constant $\lambda$. A contradiction. $\Box$ \\

\section{Cubic surface with one ordinary double point}

Now we assume that the cubic surface $X$ has one ordinary double
point $O$. If we blow up $O$, then we will get the minimal
resolution of $X$. Denote the blow up map by $\pi:\tilde{X} \to X$,
then we have $K_{\tilde{X}}=\pi^{*} K_X$. So the $\alpha_m$ and
$\alpha_{m,2}$ invariants of $X$ equal that of $\tilde{X}$. In this
section, we estimate these invariants. Since the computation on
$\tilde{X}$ is quite similar to that of
the smooth case, we shall be sketchy here.\\

\subsection{The $\alpha_m$ invariants of $\tilde{X}$}

It is well known that $\tilde{X}$ can be realized as $\C P^2$ blown
up at six ``almost general" points $p_1, \dots, p_6$. Here ``almost
general" means that three of the six points lie on a common line,
but no four of them lie on a common line and these six points are
not on a quadratic curve.(\cite{[GH]},\cite{[Dol]}) Suppose
$p_1,p_2,p_3$ lie on a common line whose strict transform is the
$(-2)$-curve $C$. We denote the exceptional divisors by $E_1, \dots,
E_6$; denote the strict transforms of the line passing through $p_i$
and $p_j$ by $L_{ij}$; and denote the strict transform of the
quadratic curve passing through $p_1, \dots, \hat{p_i},\dots, p_6$
by $F_i$. It is easy to see that there are 21 $(-1)$-curves on
$\tilde{X}$:
\begin{itemize}
  \item $E_1,E_2,E_3,E_4,E_5, E_6$;
  \item $L_{14},L_{24},L_{34}, L_{15},L_{25},L_{35},L_{16},L_{26},L_{36}, L_{45},L_{46}, L_{56}$;
  \item $F_1, F_2, F_3$.
\end{itemize}

 There are six $(-1)$-curves that intersect with the $(-2)$-curve $C$:
 $E_1,E_2,E_3$ and $L_{45},L_{46}, L_{56}$. For any of these six
 curves, there is a smooth rational curve passing though the
 intersection point of the $(-1)$-curve and $C$, and together these
 three curves constitute an anticanonical divisor of
 $\tilde{X}$. This fact in particular implies that $\alpha_m(\tilde{X})\leq
 \frac{2}{3}$. In \cite{[Chel3]}, I. Cheltsov proved $\alpha_m(\tilde{X})=\frac{2}{3}$.(See also \cite{[PaWo]} for the computation of
 $\alpha_1$.) Actually, this also follows easily from the following
 lemma:

\begin{lem}
  For any $0<\lambda\leq \frac{2}{3}$ and any nonzero holomorphic section
  $s\in H^0(X,-mK_{\tilde{X}})$, the locus of non-integrable points of
  $|s|^{-\frac{2\lambda}{m}}$ is connected. If it is not an isolated point, then it
  must be the $(-2)$-curve $C$, and in this case, $\lambda=\frac{2}{3}$ and $m$ must be even.
  Moreover, we have
  $$Z(s)=\frac{3m}{2}C+\frac{m}{2}(E_1+E_2+E_3+L_{45}+L_{46}+L_{56}).$$
\end{lem}

\noindent {\bf Proof:} It suffices to consider the case when
$\lambda=\frac{2}{3}$. As in the proof of Lemma 3.1, by Nadel's
vanishing theorem, we know that the locus of non-integrable points
of
  $|s|^{-\frac{4}{3m}}$, denoted by $LT(s)$, is a connected subset of
  $\tilde{X}$.\\

  If $LT(s)$ is not an isolated point, then we may assume that $LT(s)=\cup_i^k C_i$, where these $C_i's$
  are different irreducible curves. We can also write
  $$Z(s)=\sum_{i}^k\mu_i C_i +\Omega$$
  where $\Omega$ is an effective divisor whose support does not contain
  any of the $C_i's$. By definition of $C_i$, we have $\mu_i\geq
  \frac{3m}{2}$. It's easy to see that each $C_i$ is a smooth rational
  curve. Moreover, we have
  $$3m=-Z(s)\cdot K_{\tilde{X}}\geq -\sum_i \mu_i C_i\cdot K_{\tilde{X}}
  \geq \frac{3m}{2} \sum_i C_i\cdot (-K_{\tilde{X}}) \Rightarrow \sum_i^k C_i\cdot (-K_{\tilde{X}})\leq 2.$$
  There are three possibilities:
  \begin{enumerate}
    \item $\sum_i^k C_i\cdot (-K_{\tilde{X}})=0$;
    \item $\sum_i^k C_i\cdot (-K_{\tilde{X}})=1$;
    \item $\sum_i^k C_i\cdot (-K_{\tilde{X}})=2$.
  \end{enumerate}
  We now consider the three cases one by one.\\

  In case 1, there can be only one irreducible curve in $\sum_i
  C_i$, and it must be the $(-2)$-curve $C$. So we can write $Z(s)=\mu
  C+\Omega$, with $\mu\geq \frac{3m}{2}$. Choose any $(-1)$-curve $E$ that has
  positive intersection number with $C$. By computing the
  intersection number $E\cdot Z(s)$, it is easy to see that $\Omega$
  must contain $E$ with multiplicity at least $\frac{m}{2}$. But if this is
  true, we will have
  $$3m=\Omega \cdot (-K_{\tilde{X}})\geq 6\cdot \frac{m}{2}=3m,$$
  so we must have
  $$\Omega=\frac{m}{2}(E_1+E_2+E_3+L_{45}+L_{46}+L_{56})$$
  and $\mu=\frac{3m}{2}$. So $m$ must be even. In this case, it is easy to check that
  $$\frac{3m}{2}C+\frac{m}{2}(E_1+E_2+E_3+L_{45}+L_{46}+L_{56})\in |-mK_{\tilde{X}}|,$$
  and $LT(s)=C$.\\

  In case 2, we have either $LT(s)=C_1$ or $LT(s)=C\cup C_1$, where
  $C_1$ is a $(-1)$-curve. \\

  If $LT(s)=C_1$, we can write
  $$Z(s)=\frac{3m}{2}C_1+\Omega$$
  where $\Omega$ does not contain $C_1$. Then
  $C_1\cdot\Omega=\frac{5m}{2}$. On the other hand, $\Omega$ must
  contain any $(-1)$-curve that intersects with $C_1$, with coefficient at least $\frac{m}{2}$.
  So there are at most 5 such $(-1)$-curves. In this case, it is easy
  to see that $C_1\cdot C=1$. Thus $C\cdot\Omega=-\frac{3m}{2}$, which
  implies that $\Omega$ contains $C$ with multiplicity at least
  $\frac{3m}{4}$. But
  $$\frac{5m}{2}=C_1\cdot\Omega\geq 5\cdot \frac{m}{2}+\frac{3m}{4}> \frac{5m}{2},$$
  a contradiction.\\

  If $LT(s)=C\cup C_1$, then we can write
  $$Z(s)=\frac{3m}{2}(C+C_1)+\Omega,$$
where $\Omega$ contains neither $C$ nor $C_1$. Then
$-K_{\tilde{X}}\cdot \Omega=\frac{3m}{2}$. However, as before, we
can show that $\Omega$ must contain at least 5 $(-1)$-curves that
have positive intersection numbers with $C_1$, with multiplicity
greater than or equal to $\frac{m}{2}$. This implies that
$-K_{\tilde{X}}\cdot \Omega\geq
\frac{5m}{2}>\frac{3m}{2}$. A contradiction.\\

  In case 3, we have either $LT(s)=C_1\cup C_2$ or $LT(s)=C\cup C_1\cup
  C_2$, where $C_1$ and $C_2$ are both $(-1)$-curves. \\

  If $LT(s)=C_1\cup C_2$, then
  $$Z(s)=\frac{3m}{2}(C_1+C_2)+\Omega,$$
  where $\Omega$ contains neither $C_1$ nor $C_2$. Since $-K_{\tilde{X}}\cdot
  \Omega=0$, we have $\Omega=kC$ for some nonnegative integer $k$.
  Choose any $(-1)$-curve $L$ such that $L\cap C=\varnothing$, then
  $$m=L\cdot Z(s)=\frac{3m}{2}(L\cdot C_1+L\cdot C_2).$$
  But this is impossible, since the right hand side never equals
  $m$.\\

  If $LT(s)=C\cup C_1\cup
  C_2$, then it is easy to see that actually
  $Z(s)=\frac{3m}{2}(C+C_1+C_2)$. We can easily get a contradiction as
  above. $\Box$ \\

Since the canonical bundle of $\tilde{X}$ is nef and big, we can use
Nadel's vanishing theorem as in the proof of Theorem 3.1 to get the
following:

\begin{prop}
  For any integer $m>0$, we have $\alpha_m(\tilde{X})=\frac{2}{3}$.\\
\end{prop}

\subsection{The $\alpha_{m,2}$ invariants of $\tilde{X}$}

We have the following theorem:

\begin{thm}
  The $\alpha_{m,2}$ invariants of $\tilde{X}$ is strictly bigger
  than $\frac{2}{3}$.
\end{thm}

As in the proof of Theorem 1.2, the key to the proof of Theorem 5.1
is the following proposition:

\begin{prop}
  If $s\in H^0(\tilde{X},-mK)$ satisfies $c(s)=\frac{2}{3m}$, and $LT(s)$ is an isolated point, then $s=s_1^{\otimes
  m}$ where $s_1\in H^0(\tilde{X},-K)$ and $c(s)=\frac{2}{3}$.
\end{prop}

\noindent {\bf Proof: } Suppose $|s|^{-\frac{4}{3m}}$ is not
integrable around a point $p$. Then in view of Lemma 4.3 and Lemma
5.1, we have $$mult_p s > \frac{3m}{2}.$$ This in particular implies
that $p$ lies on $C$ or some $(-1)$-curves, for otherwise we can
find an effective divisor $D\in |-K|$ with $mult_p D=2$. By Lemma
4.1, $D\nsubseteq supp Z(s)$. So we have
$$3m=D\cdot Z(s)\geq 2 mult_p s,$$
which leads to a contradiction.\\

Moreover, if $p$ is not an Eckardt point, then as observed by
Cheltsov (Theorem 3.2 of \cite{[ChSh]}) $p$ must lie on the
$(-2)$-curve $C$. Actually, if $p \notin C$, then we can repeat the
proof of Theorem 4.1 to show that $p$ is an Eckardt point.\\

Now we prove that $p$ also lies on a $(-1)$-curve. Suppose not, then
it is easy to see that there is an irreducible curve $D$ with
$mult_p D=2$ and $D+C \in |-K|$. By Lemma 4.1, $D\nsubseteq supp
Z(s)$. So we have
$$3m=D\cdot Z(s)\geq 2 mult_p s>3m.$$
A contradiction.\\

Denote the $(-1)$-curve through $p$  by $L$, then there is a unique
irreducible curve $D$ such that $p\in D$ and $C+L+D\in |-K|$. Write
$Z(s)$ as $Z(s)=\mu C+ \nu L + \Omega$, where  $\Omega$ is an
effective divisor whose support contains neither $C$ nor $L$. Asumme
$D$ is not contained in the support of $\Omega$. We shall
derive a contradiction from this assumption.\\

First, by Lemma 4.3 and Lemma 5.1, we know that
$$
  mult_p s > \frac{3m}{2}.
$$
By our assumption, we also have $D \cdot \Omega>0$, $C\cdot
\Omega>0$ and $L\cdot \Omega>0$. A careful analysis of these three
inequalities leads to the following results:
$$mult_p s\leq 2m,\quad mult_p \Omega \leq \frac{5m}{6}, \quad \mu >\frac{m}{2}, \quad \nu >\frac{m}{4}.$$
Choose a suitable blow down map $\pi: \tilde{X}\to \C P^2$ such that
$deg \pi (C)=1=deg \pi(L)=deg \pi(D)$. Then five out of the six
blowing up centers lie on $\pi(C)\cup\pi(L)$, with the other one,
denoted by $q$, lying on $\pi(D)$. By Lemma 4.1 and our assumption,
$\pi(\Omega)$ can not contain any line through $q$. By computing
intersection numbers of $\pi(\Omega)$ with lines through $q$, we can
get:
$$\mu\leq m, \quad \nu\leq m.$$

We now further blow up $\tilde{X}$ at $p$ with exceptional divisor
$E$. The blowing up map is $f: X_1 \to \tilde{X}$. Then
$$f^{\ast}\Big{(} K_{\tilde{X}}+ \frac{2}{3m} Z(s)\Big{)}=
K_{X_1}+\frac{2}{3m}\Big{(}\mu \bar{C}+\nu\bar{L}+\bar{\Omega}
\Big{)}+(\frac{2}{3m}mult_p s-1)E,$$ where $\bar{F}$ denotes the
strict transform of $F$ for any divisor $F$. Since $\frac{2}{3m}
Z(s)$ is not log terminal at $p$, there exists a point $Q$ on $E$
such that $\frac{2}{3m}\Big{(}\mu \bar{C}+\nu\bar{L}+\bar{\Omega}
\Big{)}+(\frac{2}{3m}mult_p s-1)E$ is not log terminal at $Q$.\\

If $Q\notin \bar{C}$ and $Q\notin \bar{L}$, then
$$\frac{2}{3m} mult_Q \bar{\Omega}+(\frac{2}{3m}mult_p s-1)\geq 1,$$
hence $mult_Q \bar{\Omega} \geq 3m- mult_p s \geq m.$ But this is
impossible in view of the fact $mult_p \Omega \leq \frac{5m}{6}$.
Now suppose $Q \in \bar{L}$. Then by Fubini's theorem on repeat
integration, we have
$$\bar{L} \cdot \Big{(} \frac{2}{3m}\bar{\Omega}+(\frac{2}{3m}mult_p s-1)E \Big{)}\geq 1,$$
which implies that $\nu\geq m$. So we have $\nu=m$ and $mult_p
\Omega \leq m-\mu< \frac{m}{2}$. However, by H\"older inequality,
this implies $c_p(s)>\frac{2}{3m}$.\\

So there is only one possibility, namely, $Q \in \bar{C}$. As in the
above discussions, we have
$$\bar{C}\cdot \Big{(} \frac{2}{3m}\bar{\Omega}+(\frac{2}{3m}mult_p s-1)E \Big{)}\geq 1.$$
This implies $\mu\geq m$, hence $\mu=m$. In this case, we have
$\nu=L\cdot \Omega \geq mult_p \Omega$ and $mult_p \Omega+\nu=mult_p
s-m\leq m$. Combined with H\"older inequality, these inequalities
imply that $\nu=mult_p \Omega =\frac{m}{2}$. This at once leads to a
contradiction when $m$ is odd. When $m$ is even, this is also
impossible by the following lemma:

\begin{lem}
  Let $f,h$ be germs of holomorphic functions at $0\in \C^2$,
  with $f(z_1,z_2)=z_1^{2k} z_2^k h$, and $mult_0 h=k$. Suppose
  $z_1 \nmid h$ and $z_2 \nmid h$ in $\mathcal {O}_{\C^2,0}$, then
  $c_0(f)>\frac{1}{3k}$.
\end{lem}
This lemma can be proved, for example, by induction on the number of
blowing ups to resolve the singularity of $\{f=0\}$. The detail is
left as an exercise for the reader. $\Box$ \\

\subsection{K\"ahler-Einstein metrics on singular cubic surfaces}

Even though we have
$$\frac{1}{\alpha_m(X)}+\frac{1}{\alpha_{m,2}(X)}>3$$ for $X$, we can
not apply Theorem 1.2 directly to show the existence of
K\"ahler-Einstein metrics on $X$ due to the presence of
singularities. We now explain these difficulties.\\

Recall that in Tian's proof of the Calabi's conjecture in dimension
2 \cite{[T3]}, he first found one surface in the moduli space that
admits a K\"ahler-Einstein metric (this was done in an earlier paper
with Yau \cite{[TY]}), then he used the continuity method and solved
the Monge-Amp\`ere equations along a regular family of complex
surfaces. The key point is the $C^0$-estimate, for the higher order
estimates differ little from Yau's paper \cite{[Yau]}. To prove the
$C^0$-estimate, there are three main ingredients: the ``partial
$C^0$-estimate", the $\alpha_m, \alpha_{m,2}$-invariants estimates
and the continuity of rational integrals in dimension 2. The
``partial $C^0$-estimate" is used to reduce the $C^0$-estimate to
the uniform estimate of certain rational integrals. Then we can
combine the $\alpha$-invariants
estimates and the continuity of rational integrals to get the uniform $C^0$-estimate.\\

In our case, one may do everything in the category of orbifolds. We
can start with a family of cubic surfaces, each has one ordinary
double points, together with a family of smooth varying orbifold
K\"ahler metrics. The analysis here is almost identical to that of
the smooth case considered in \cite{[T3]}, but the problem is that
we do not {\em a priori} have the existence of one K\"ahler-Einstein
orbifold in the moduli space of cubic surfaces with one ordinary
double point in view of our computations in Lemma 5.1. So instead,
we choose a family of smooth cubic surfaces degenerating to the
singular surface $X$. We solve the Monge-Amp\`ere equations along
this family. However, the $C^2$-estimate does not follow easily from
the $C^0$-estimate in our case, since we do not have a uniform lower
bound for the bisectional curvature. Actually, we can expect a
``partial $C^2$-estimate", i.e., a $C^2$-estimate outside a
subvariety. At this point, the technique of \cite{[ST]} and
\cite{[To]} should be helpful. Also, we need a generalized
continuity theorem to insure the continuity of rational integrals
when the integration domain degenerates to a domain with ``mild"
singularities.\\

For general normal cubic surfaces, a theorem of Ding and Tian
\cite{[DiTi]} claims that if the surface has a K\"ahler-Einstein
metric, then it must be semistable. That is to say, the surface can
not have singularities other than $A_1$ and $A_2$ types. It is
expected that the main theorem of \cite{[T7]} still holds in the
orbifold case, then we can easily show that a semistable cubic
surface with a $A_2$ singular point can not have K\"ahler-Einstein
metrics except that it has three $A_2$ singular points. Note that
every cubic surface with three $A_2$ singular points is projectively
equivalent to the surface defined by $$z_0^3+z_1z_2z_3=0.$$ It is
the quotient of $\C P^2$ by the cyclic group $\Gamma_3$, hence it
always has an orbifold K\"ahler-Einstein metric (See
\cite{[DiTi]}).\footnote{In \cite{[Chel3]}, Cheltsov also proved the
existence of K\"ahler-Einstein metric by computing the
$\alpha_G$-invariant .}\\

A cubic surface with only $A_1$ singularities (i.e. ordinary double
points) can be easily classified as did in \cite{[BrWa]},
\cite{[GH]} and \cite{[Dol]}. The number of singular points could be
1, 2, 3 or 4. Any cubic surface with 4 ordinary double points is
projectively equivalent to the Segre cubic surface defined by the
equation
$$z_0z_1z_2+z_0z_1z_3+z_0z_2z_3+z_1z_2z_3=0.$$ In \cite{[Chel3]}
Cheltsov proved the existence of orbifold K\"ahler-Einstein metric
on this surface by showing that the
$\alpha_G$-invariant is bigger than 2/3. \\

If the cubic surface $X$ has 2 or 3 ordinary double points, then
$\alpha_m(X)=1/2$ according to Cheltsov \cite{[Chel3]}. One can
easily show that $\alpha_{m,2}(X)\leq 1$. So even if we can overcome
all the analytic difficulties mentioned above, we still can not use
Tian's criterion. In this case, one needs a generalized form of
Theorem 1.2, also making use of $\alpha_{m,3}$. Note that for a
cubic surface with 1 or 2 ordinary double points, the
$\alpha_G$-invariants are equal to the corresponding
$\alpha$-invariants, due to the existence of certain $G$-invariant
anticanonical divisors.\\

All these problems will be discussed in details in another paper.\\

\appendix

\section{Appendix}

In this appendix, we give a proof of Theorem 2.2. This part is taken
from \cite{[Shi]}.\\

\noindent {\bf Proof of Theorem 2.2:} Since we always have
$$\alpha(X)\leq \inf_m \alpha_m(X),$$ it's clear that we need only to
 prove: $\forall c>\alpha(X)$, we can find some $m\gg 1$ s.t. $c>\alpha_m(X)$.\\

Now by definition, we can find a sequence of $\varphi_i \in P(X,g)$,
s.t. $$\int_X e^{-c\varphi_i}dV_g \to +\infty,\quad as\
i\to+\infty$$ To make things simple, we assume further that $\exists
\varepsilon_0>0$, s.t.
$$\int_X e^{-(c-\varepsilon_0)\varphi_i}dV_g \to +\infty,\quad as\ i\to+\infty$$ Set $\bar{c}=c-\frac{\varepsilon_0}{2}$.
After passing to subsequence, we may assume without loss of
generality that
 $\varphi_i$ converges to a $\varphi$ in $L^1_{loc}$. Then $\omega_g+\frac{\sqrt{-1}}{2\pi} \partial \bar{\partial} \varphi \geq 0$.\\

Define $$c(\varphi):= \sup \{\alpha>0 \big{|}\  \int_X
e^{-\alpha\varphi} dV_g <+\infty \} $$

{\claim $\bar{c}>c(\varphi)$}

The reason is simple: If $\bar{c}\leq c(\varphi)$, then $\forall
\varepsilon>0,\ \bar{c}-\varepsilon<c(\varphi)$. Then by Theorem
2.1, we know that $$\int_X e^{-(\bar{c}-\varepsilon)\varphi_i}dV_g
\to \int_X e^{-(\bar{c}-\varepsilon)\varphi}dV_g< +\infty,\quad as\
i\to+\infty$$
A contradiction.\\

{\claim $\exists k_0>0,\ s.t.\ \forall m>k_0$, there exists a global
non-zero holomorphic section $s$ of $K^{-m}$, such that $$\int_X
|s|^2_{h^m} e^{-(m-k_0)\varphi}dV_g < +\infty$$ }

We will first prove the theorem assuming this.\\

$\forall p>1, q>1$ with $\frac{1}{p}+\frac{1}{q}=1$, we have
{\beqs \int_X e^{-(m-k_0)p^{-1}\varphi}dV_g&=&\int_X \big{(}|s|^2_{h^m} e^{-(m-k_0)\varphi}\big{)}^{\frac{1}{p}}\big{(}|s|^2_{h^m}\big{)}^{\frac{1}{q}-1}dV_g\\
& \leq & \Huge{(}\int_X |s|^2_{h^m} e^{-(m-k_0)\varphi}dV_g
\Huge{)}^{\frac{1}{p}}\ \Huge{(}\int_X |s|^{-2(q-1)}_{h^m} dV_g
\Huge{)}^{\frac{1}{q}} \eeqs}

If we let $q-1=\frac{\bar{c}}{m}$, then
$\frac{m-k_0}{p}=\frac{\bar{c}(m-k_0)}{m+\bar{c}}<\bar{c}$. When $m$
is big enough, $\frac{m-k_0}{p}$ is very close to $\bar{c}$, hence
$\frac{m-k_0}{p}> c(\varphi)$. So we have
$$\int_X |s|^{-\frac{2\bar{c}}{m}}_{h^m} dV_g =+\infty$$
That is to say $c>\bar{c}\geq \alpha_m (X)$. So we have
$\alpha(X)=\inf_m \alpha_{m}(X)$.\\

It also can be seen from the above proof that for any sequence
$m_k\nearrow +\infty$, $$\alpha(X)=\inf_k \alpha_{m_k}(X).$$ From
this, we have $$\liminf_{m\to +\infty}\alpha_m(X)=\alpha(X).$$ If
$$\limsup_{m\to +\infty}\alpha_m(X)=c>\alpha(X),$$ then we have a
sequence $m_k\nearrow +\infty$, s.t. $\alpha_{m_k}(X)\to c$. Then
when $k$ is sufficiently large, say, larger than some fixed $k_0$,
we
 have $$\alpha_{m_k}(X)>\frac{c+\alpha(X)}{2}.$$But this is absurd,
 since we still have $$\alpha(X)=\inf_{k>k_0} \alpha_{m_k}(X).$$ So we must have
$\alpha(X)=\lim_{m\to +\infty} \alpha_{m}(X)$. $\Box$\\

\noindent {\bf Proof of Claim 2:} The proof is standard. It is
motivated by a theorem of Demailly (\cite{[De]} P110-111). For
completeness and the
reader's convenience, we include a proof here following \cite{[De]}. \\

Fix $x\in X$, s.t. $\varphi(x)\neq -\infty$. Choose a pseudoconvex
coordinate neighborhood $\Omega$ of $x$, such that $K^{-m}$ is
trivial over $\Omega$. Then by the Ohsawa-Takegoshi
theorem\footnote{In fact, H\"ormander's $L^2$ theory is enough. This
result is known as the H\"ormander-Bombieri-Skoda theorem in the
literature.}, there exists a holomorphic function $g$ on $\Omega$
with $g(x)=1$ and
$$\int_\Omega |g(z)|^2e^{(m-k_0)\varphi(z)}d\lambda(z)<C_1$$($k_0$
to be fixed later). Choose a $C^\infty$ cut-off function $\chi : \R
\to [0,1]$ s.t. $\chi(t)|_{t\leq 1/2}\equiv 1,\ \chi(t)|_{t\geq
1}\equiv 0$. Using the trivialization of $K^{-m}$, we may view
$\chi(\frac{z-x}{r})g(z)
$ as a smooth section of $K^{-m}$. We denote it by $\sigma$ to avoid confusion. Put $v:=\bar{\partial}\sigma$. We want to solve the equation $\bar{\partial}u=v$.\\

To make sure the section we shall construct is non-trivial, we
introduce a new weight $e^{-2n\rho_x(z)}$ where
$\rho_x(z)=\chi(\frac{z-x}{r})\log|z-x|$. Choose $k_0\in \N$ s.t.
$k_0 \omega_g + n\frac{\sqrt{-1}}{\pi}\partial
\bar{\partial}\rho_x\geq 0$. Then we choose a singular hermitian
metric on $K^{-m}$ to be $h^me^{-(m-k_0)\varphi}$. Since
$$n\frac{\sqrt{-1}}{\pi}\partial \bar{\partial}\rho_x+Ric(h)+Ric(h^me^{-(m-k_0)\varphi})\geq \omega_g,$$
by H\"ormander's $L^2$ existence theorem, we get a smooth section
$u$ with $\bar{\partial}u=v$ and
$$\int_X |u|^2_{h^m}e^{-2n\rho_x-(m-k_0)\varphi}dV_g\leq C_2 \int_{\frac{1}{2}<|z-x|<r}|g|^2e^{-2n\rho_x-(m-k_0)\varphi}dV_g\leq C_3$$
Since $\rho_x$ has logarithmic singularity at $x$, we have $u(x)=0$.
Define $s=\sigma-u$, then $s$ is a non-zero holomorphic section of
$K^{-m}$ with $$\int_X |s|^2_{h^m} e^{-(m-k_0)\varphi}dV_g <
+\infty.$$
 $\Box$\\

It's easy to see that our proof still works in the
orbifold case.\\

\begin{flushleft}
School of Mathematical Sciences, Peking University,\\ Beijing, 100871, P. R. China.\\
Email: ylshi@math.pku.edu.cn
\end{flushleft}

\end{document}